\documentclass[11pt,reqno]{amsart}
\usepackage{amssymb}
\usepackage{amsfonts}
\usepackage{amsthm}
\numberwithin{equation}{section}
\theoremstyle{plain}
\newtheorem{theorem}{Theorem}[section]
\newtheorem{proposition}[theorem]{Proposition}
\newtheorem{corollary}[theorem]{Corollary}
\newtheorem{lemma}[theorem]{Lemma}
\theoremstyle{definition}
\newtheorem{definition}[theorem]{Definition}
\newtheorem*{remark}{Remark}
\newtheorem*{example}{Example}
\newcommand{\refE}[1]{(\ref{E:#1})}
\newcommand{\refS}[1]{Section~\ref{S:#1}}
\newcommand{\refSS}[1]{Section~\ref{SS:#1}}
\newcommand{\refT}[1]{Theorem~\ref{T:#1}}

\newcommand{\refP}[1]{Proposition~\ref{P:#1}}
\newcommand{\refD}[1]{Definition~\ref{D:#1}}
\newcommand{\refC}[1]{Corollary~\ref{C:#1}}
\newcommand{\refL}[1]{Lemma~\ref{L:#1}}
\newcommand{\R}{\ensuremath{\mathbb{R}}}
\newcommand{\C}{\ensuremath{\mathbb{C}}}
\newcommand{\N}{\ensuremath{\mathbb{N}}}

\renewcommand{\P}{\ensuremath{\mathbb{P}}}
\newcommand{\Z}{\ensuremath{\mathbb{Z}}}
\newcommand{\K}{\ensuremath{\mathbb{K}}}

\newcommand{\cintl}[1]{\frac 1{24\pi\mathrm{i}}\int_{#1 }}
\newcommand{\g}{\ensuremath{\gamma}}

\newcommand{\A}{\mathcal{A}}
\newcommand{\W}{\mathcal{W}}
\newcommand{\Wh}{\widehat{\mathcal{W}}}
\newcommand{\Vir}{\mathcal{V}}
\renewcommand{\L}{\mathcal{L}}

\newcommand{\Lh}{\widehat{\mathcal{L}}}

\newcommand{\fpz}{\frac {d }{dz}}
\newcommand{\pfz}[1]{\frac {d#1}{dz}}

\renewcommand{\H}{\mathrm{H}}
\renewcommand{\a}{\ensuremath{\alpha}}
\renewcommand{\b}{\ensuremath{\beta}}
\newcommand{\Bh}{\widehat{B}}
\newcommand{\E}{\mathcal{E}}
\newcommand{\Dex}{D_1^*\cup D_{-1/2}^*\cup D_{-2}^*}
\newcommand{\eee}{(e_1-e_2)(e_1-e_3)}
\newcommand{\ddX}{\frac {d}{dX}}

\newcommand{\tensor}{\otimes}   
\newcommand{\htensor}{\widehat{\otimes}}   
\newcommand{\gl}{\mathfrak{gl}} 
\begin{document}
\vspace*{-1cm}
\hbox{ }
{{\hspace*{\fill} Mannheimer Manuskripte 266}}

{{\hspace*{\fill} math.QA/0206114}}

{{\hspace*{\fill} MPI - Bonn, Nr. 2002-53 }}

\vspace*{2cm}

\title{Global Deformations of the Witt Algebra of
Krichever-Novikov Type}
\thanks{The research of the 
first author was partially supported by the grants OTKA
T030823, T29525. The last version of the work was completed during her
stay at the Max-Planck-Institut f\"ur Mathematik Bonn.}

\author[A. Fialowski]{Alice Fialowski}
\address[Alice Fialowski]{Department of Applied Analysis,
E\"otv\"os Lor\'and University, P\'azm\'any P\'eter s\'et\'any 1,
H-1117 Budapest, Hungary}
\email{fialowsk@cs.elte.hu}
\author[M. Schlichenmaier]{Martin Schlichenmaier}
\address[Martin Schlichenmaier]{Department of Mathematics and
  Computer Science, University of Mannheim, A5,
         D-68131 Mannheim,
          Germany}
\email{schlichenmaier@math.uni-mannheim.de}
\begin{abstract}
By considering non-trivial global deformations of the Witt (and the
Virasoro) algebra given by geometric constructions it is shown 
that, despite their infinitesimal and formal rigidity, they
are globally not rigid.
This shows the need of a clear indication of the type of deformations
considered.
The families appearing are constructed as families of 
algebras of Krichever-Novikov
type. 
They show up in 
a natural way in the global operator approach 
to the quantization of two-dimensional conformal field theory.
In addition, a proof of the infinitesimal and formal rigidity 
of the Witt algebra is presented.
\end{abstract}
\date{June 4, 2002}
\maketitle
\vskip 0.5cm
\hskip 2cm
\begin{minipage}{0.6\linewidth}
{\it We dedicate this article to the memory of our
good friend  Peter Slodowy who passed away in 2002} 
\end{minipage}

\vskip 1.0cm
\section{Introduction}\label{S:intro}
Deformations of mathematical structures are important in most
part of mathematics and its applications. Deformation is one of the
tools used to study a specific object, by
deforming it into some families of "similar" structure objects.  
This way we get a richer picture about the original object itself.

But there is also another question approached via deformations.
Roughly speaking, it is the question,
can we equip the set of mathematical structures under consideration
(maybe up to certain equivalences) with the structure of
a topological or even geometric space.
In other words, does there exists a moduli space for these
structures.
If so, then for a fixed object the deformations of this object should
reflect the local structure of the moduli space at 
the point corresponding to this object.

In this respect, 
a clear success story  is the classification of complex
analytic structures on a fixed topological manifold.
Also in algebraic geometry one has well-developed results in this 
direction.
One of these results is that the local situation at a point $[C]$
of the moduli space is completely governed by the cohomological 
properties of the geometric object $C$.
As typical example recall that for the moduli space $\mathcal{M}_g$ of 
smooth projective curves of genus $g$ over $\C$ (or equivalently, compact
Riemann surfaces of genus $g$)
the tangent space $T_{[C]} \mathcal{M}_g$ can be naturally 
identified with 
$\H^1(C,T_C)$, where $T_C$ is the sheaf of holomorphic vector fields
over $C$.
This extends to higher dimension. In particular, it turns out that
for compact complex manifolds $M$,
the condition $\H^1(M,T_M)$ implies that $M$ is rigid, \cite[Thm. 4.4]{Kod}.
Rigidity means  that
any differentiable family $\pi:\mathcal{M}\to B\subseteq \R$, 
$0\in B$ which contains $M$ as the special member $M_0:=\pi^{-1}(0)$
 is trivial in a neighbourhood of $0$, i.e. for $t$ small enough
$M_t:=\pi^{-1}(t)\cong M$.
Even more generally, for $M$ a compact complex manifold
and $\H^1(M,T_M)\ne \{0\}$ there exists a versal family 
which can be realized locally as a family  over a certain subspace
of    $\H^1(M,T_M)$ such that every appearing deformation family
is ``contained'' in this versal family
(see also \cite{Palm} for definitions,  results, and further references).

These positive results lead to the impression that the
vanishing of the relevant cohomology spaces will imply rigidity 
with respect to deformations also in the case of other
structures.

The goal of this article is to shed some light on this in the context
of deformations of infinite-dimensional Lie algebras.
We consider the case of the Witt algebra $\W$ (see its definition 
further down).
The  cohomology 
``responsible'' for deformations is $\H^2(\W,\W)$.
It is known that $\H^2(\W,\W)=\{0\}$ (see \refS{rigid}).
Hence, guided by the experience in the theory of deformations
of complex manifolds, one might think that $\W$ is rigid in the sense
that
all families containing $\W$ as a special element will be 
trivial.
But this is not the case as we will show.
Certain natural families
of Krichever-Novikov type algebras of geometric origin (see \refS{kn}
for their definition)
will appear which  contain  $\W$ as special element.
But none of the other elements will be isomorphic to $\W$.
In fact, from  $\H^2(\W,\W)=\{0\}$
it follows that the Witt algebra is infinitesimally and formally rigid.
But this condition does not imply that there are no non-trivial 
global deformation families.
The main point to learn is that 
it is necessary to distinguish clearly the formal and the global
deformation situation.
The formal rigidity of the Witt algebra indeed follows from 
$\H^2(\W,\W)=\{0\}$,  but no statement like that about 
global deformations.

How intricate the situation is, can be seen from the fact that
for the subalgebra $L_1$ of $\W$,
consisting of those vector fields vanishing at least of order two at
zero, 
there exists a versal formal family consisting of three different families
parameterized over a collection of three curves in 
$\H^2(L_1,L_1)$.
Suitably adjusted each family  corresponds 
to a scalar multiple of the same  
cohomology class $\omega$  \cite{Fiaproc},\cite{Fiajmp},\cite{FiaFuMini}
which gives their infinitesimal deformations.
It turns out that $\H^2(L_1,L_1)$ is 3-dimensional, but only the
infinitesimal deformations corresponding to scalar multiples of   $\omega$
can be extended to formal deformations.

The results of this article will show that the theory of
deformations of infinite-dimensional Lie algebras is still not
in satisfactory shape. 
Hopefully, the appearing features will be of help in a further
understanding.

Clearly, what will be done here, can also be done in the case of
associative algebras.
In particular, there will be global deformations of the
associative algebra of Laurent polynomials of
Krichever-Novikov type obtained by the same
process as the one presented here.

First, let us introduce the basic definitions.
Consider the complexification of the Lie algebra of polynomial vector
fields on the circle with generators 
$$
l_n:= \exp (\mathrm{i} n\varphi)\frac{d}{d\varphi}, \qquad n \in \Z,
$$
where $\varphi$ is the angular parameter.  The bracket operation in this
Lie algebra is
$$
[l_n,l_m] = (m-n)l_{n+m}.
$$
We call it the \textit{Witt algebra} and denote it by $\W$.
Equivalently, the Witt algebra can be described as the Lie algebra of
meromorphic vector fields on the Riemann sphere $\P^1(\C)$ which are
holomorphic outside the points $0$ and $\infty$.
In this presentation $l_n=z^{n+1}\fpz$, where $z$ is the quasi-global
coordinate on  $\P^1(\C$).

 The   Lie algebra $\W$ is infinite
dimensional and graded with the standard grading 
$\deg l_n =n$.
By taking formal vector fields with the projective limit topology we get
the completed topological Witt algebra $\Wh$. Throughout this paper we
will consider its everywhere dense subalgebra $\W$. 

It is well-known that $\W$ has a unique nontrivial one-dimensional
central extension, the \textit{Virasoro algebra} $\Vir$. 
It is generated by $l_n\ (n \in \Z)$ and the
central element $c$, and its bracket operation is defined by 
\begin{equation}\label{E:Vir}
[l_n,l_m] = (m-n)l_{n+m} + \tfrac{1}{12}(m^3-m)\delta_{n,-m}\,c,
\qquad
[l_n,  c] =0.
\end{equation}
\medskip

In \refS{theory} we recall the different concepts of deformation.
There is a lot of confusion in the literature in the notion of a
deformation. Several different (inequivalent) approaches exist. One of
the aims of this paper is to clarify the difference between
deformations of geometric origin and so called formal deformations. Formal
deformation theory has the advantage of using cohomology. It is also
complete in the sense that under some natural cohomology assumptions
there exists a versal formal deformation, which induces all other
deformations.

In this context see \refT{formver} and \refC{formrig}
which are quoted from \cite{Fiaproc},\cite{FiaFuMini}.
Formal deformations are deformations with a complete
local algebra as base.
A deformation with a commutative (non-local) algebra base gives a much
richer picture of deformation families, depending on the augmentation of
the base algebra. If we identify the base of deformation -- which is a
commutative algebra of functions -- with a smooth manifold, an augmentation
corresponds to
choosing a point on the manifold. So choosing different points should in 
general lead to different deformation situations.  
As already indicated above, in infinite dimension, there is no tight relation
between global and formal deformations, as we will show in this paper.

\medskip
In \refS{rigid}, we supply a detailed proof of the 
infinitesimal and formal rigidity of the Witt algebra $\W$, by 
showing that $\H^2(\W,\W)=0$.

\medskip
In \refS{kn} we introduce and recall the necessary properties of the
Krichever-Novikov type vector field algebras.
They are generalizations of the Witt algebra (in its presentation
as vector fields on $\P^1(\C)$) to higher genus smooth projective curves.

\medskip
In \refS{elliptic} 
we construct  global deformations of the Witt algebra 
by considering  certain families of algebras for
the genus one case (i.e. the elliptic curve case) 
and let the elliptic curve degenerate to
a singular cubic.
The two points, where poles are allowed, are the zero element of the
elliptic curve (with respect to its additive structure) and a
2-torsion point.
In this way we obtain families parameterized over the affine line with
the peculiar behaviour that every family is a global deformation of the
Witt algebra, i.e. $\W$ is a special member, whereas all other members
are mutually 
isomorphic but not isomorphic to $\W$, see \refT{families}. Globally
these families are non-trivial, but infinitesimally and formally they
are trivial. 
The construction can be extended to the centrally extended 
algebras, yielding a global deformation of the Virasoro algebra.
Finally, we consider the subalgebra $L_1$ of $\W$ corresponding to
the vector fields vanishing  at least with order two at $0$. This
algebra is formally not rigid, and its formal versal deformation has
been determined \cite{Fiaproc, FiaFusd}. We identify one of the
appearing three families in our geometric context.

\medskip 
The results obtained do not have only relevance to the
deformation theory of algebras but also to the theory of
two-dimensional conformal fields and their quantization.
It is well-known that the Witt algebra, the Virasoro algebra, and their representations 
are of fundamental importance for the local description of 
conformal field theory on the Riemann sphere
(i.e. for genus zero), see \cite{BP}.
Krichever and Novikov \cite{KNFac} proposed in the case of
higher genus Riemann surfaces the use of global operator fields 
which are given with the help of  the Lie algebra 
of vector fields of Krichever-Novikov type,
certain related algebras, 
and their representations (\refS{kn}).

In the process of quantization of the conformal fields one has
to consider families of algebras and representations over the
moduli space of compact Riemann surfaces
(or equivalently, of  smooth projective curves over $\C$) of genus $g$ with
$N$ marked points.
See \cite{SchlShWz} for a global operator version,
and \cite{TUY} for a sheaf version.
In passing to the boundary of the moduli space one obtains 
the limit objects which are defined over the normalization
of curves of lower genus.
Assuming good behaviour of the 
examined objects under deformation also in the
limit (e.g. ``factorization''), the degeneration
is an important technique to obtain  via induction 
from results for lower genus  results for all 
genera.
See Tsuchiya, Ueno, and 
Yamada's proof of the Verlinde formula \cite{TUY}
as an application of this principle.

By a maximal degeneration a collection of $\P^1(\C)$'s will 
appear.
For the vector field algebras (with or without central extension)
one obtains families of algebras which are related to the 
Witt or Virasoro algebra
or certain subalgebras  respectively.
Indeed, the examples considered in this article are exactly of this type.
They appear as families which are naturally defined over the
moduli space of complex one-dimensional tori (i.e. of elliptic curves)
with two marked points.
The full geometric picture behind it was discussed in \cite{SchlDeg}.
In special cases the Witt and Virasoro algebra appear as 
degenerations of the Krichever-Novikov algebras. 
Considered from the opposite point of view,
in the sense of this article,  
the Krichever-Novikov algebras
are global deformations of the Witt and Virasoro algebra.
Nevertheless, as we show here, the structure of 
these algebras are not determined 
by the Witt algebra, despite the formal rigidity of
the latter.
\section{Deformations and formal deformations}\label{S:theory}
\subsection{Intuitively}
Let us start with the intuitive definition. Let $\L$ be a Lie algebra
with Lie bracket $\mu_0$ over a field $\K$. A {\it deformation} of
$\L$ is a one-parameter family $\L_t$ of Lie algebras with the bracket
\begin{equation}\label{E:intdef}
\mu_t = \mu_0 + t\phi_1 + t^2\phi_2 + ...
\end{equation}
where the $\phi_i$ are $\L$-valued 2-cochains, i.e. elements of
$\mathrm{Hom}_{\K}(\bigwedge^2\L,\L)=C^2(\L,\L)$, and $\L_t$ is a Lie
algebra for each $t\in \K$. (see \cite {Fiasbor, Ger}). Two deformations
$\L_t$ and $\L'_t$ are equivalent if there exists a linear automorphism
$\hat\psi_t = \text{id} + \psi_1t + \psi_2t^2 + ...$ of $\L$ where $\psi_i$
are linear maps over $\K$, i.e. elements of $C^1(\L,\L)$, such that
\begin{equation}\label{E:intequ}
\mu'_t(x,y) = \hat\psi_t^{-1}(\mu_t(\hat\psi_t(x),\hat\psi_t(y)).
\end{equation}

The Jacobi identity for the algebra $\L_t$ implies that the 2-cochain
$\phi_1$ is indeed a cocycle, i.e. it fulfills
$d_2\phi_1=0$ with respect to the Lie algebra cochain complex
of $\L$ with values in $\L$ (see \cite{Fuks} for the definitions).
If $\phi_1$ vanishes identically, the first nonvanishing $\phi_i$ will
be a cocycle.
If $\mu'_t$ is an equivalent deformation (with cochains $\phi'_i$) then
\begin{equation}
\phi'_1-\phi_1=d_1\psi_1.
\end{equation}
Hence, every equivalence class 
of deformations defines uniquely an element of $\H^2(\L,\L)$.
This class is called the {\it differential} of the deformation.
The differential of a family which is
equivalent to a trivial family will be the zero cohomology class.


\subsection{Global deformations} 
Consider now a deformation $\L_t$ not
as a family of Lie algebras, but as a Lie algebra over the 
algebra $\K[[t]]$. The natural generalization is to allow more parameters,
or to take in general a commutative algebra $A$ 
over $\K$ with identity as base
of a deformation. 

In the following we will assume that $A$ is a
commutative algebra over  the field $\K$ of  characteristic zero
which  admits an augmentation 
 $\epsilon: A \to \K$.
This says that $\epsilon$ is a $\K$-algebra homomorphism, e.g.
$\epsilon(1_A)=1$.
The ideal $m_\epsilon:= \text{Ker\,} \epsilon$ is 
a maximal ideal of $A$. 
Vice versa, given a maximal ideal $m$ of $A$ with 
$A/m\cong\K$, then the natural factorization map defines an
augmentation. 

In case that $A$ is a finitely generated $\K$-algebra over
an algebraically closed field $\K$ then $A/m\cong\K$
is true for every maximal ideal $m$. Hence in this case 
every such $A$ admits at least one augmentation and
all maximal ideals are coming from augmentations.

Let us consider a Lie algebra $\L$ over the field $\K$, 
 $\epsilon$ a fixed augmentation of $A$, 
and
$m=\text{Ker\,} \epsilon$ the associated maximal ideal. 
\begin{definition}
\label{D:glob}
 A {\it global deformation} $\lambda$ of $\L$ with the
base $(A,m)$ or simply with the base $A$, is a Lie $A$-algebra structure
on the tensor product $A\tensor_\K\L$ with bracket $[.,.]_\lambda$ such
that
\begin{equation}
\epsilon \textstyle{\tensor} \text{id} : A \textstyle{\tensor} \L \to
 \K \textstyle{\tensor} \L = \L
\end{equation}
is a Lie algebra homomorphism
(see \cite{FiaFuMini}).
Specifically, it means that for all $a,b \in A$ and $x,y \in \L$,
\smallskip

\quad (1) $[a \tensor x, b \tensor y]_\lambda = (ab \tensor \text{id\,})[1 \tensor x,
1 \tensor y]_\lambda$,
\smallskip

\quad (2) $[.,.]_\lambda$ is skew-symmetric and satisfies the Jacobi 
identity,
\smallskip

\quad (3) $\epsilon \tensor \text{id\,}([1 \tensor x, 1 \tensor y]_\lambda) = 1
\tensor [x,y]$.
\end{definition}

\noindent 
By Condition (1) to describe a deformation 
it is enough to give  the elements $[1 \tensor x, 1
\tensor y]_\lambda$ for all $x,y \in \L$. 
By condition (3) follows  that for them the Lie product has the form
\begin{equation}\label{E:defex}
[1 \textstyle{\tensor} x, 1 \textstyle{\tensor} y]_\lambda = 1
\textstyle{\tensor} [x,y] + \sum_i a_i \textstyle{\tensor} z_i,
\end{equation}
with $a_i \in m$, $z_i \in \L$.

A deformation is called {\it trivial} if $A\tensor_{\K} \L$
carries the trivially extended Lie structure, i.e. \refE{defex} reads as 
$[1\tensor x,1\tensor y]_\lambda=1\tensor[x,y]$.
Two deformations of a Lie algebra $\L$ with the same base $A$ are
called {\it equivalent} if there exists a Lie algebra isomorphism
between the two copies of $A \tensor \L$ with the two Lie algebra
structures, compatible with $\epsilon \tensor \text{id}$.

We say that  a deformation
is {\it local} if $A$ is a local $\K$-algebra 
with unique maximal ideal $m_A$. 
By assumption $m_A=\text{Ker\,} \epsilon$
and $A/m_A \cong \K$.
In case that in addition  ${m_A}^2 = 0$, the deformation
is called  {\it infinitesimal}.
\begin{example}
If $A = \K[t]$, then this is the same as an algebraic 
1-parameter deformation of $\L$.
In this case we sometimes use simply the expression ``deformation
over the affine line.''
This can be extended to the case where $A$ is the algebra of regular
functions on an affine variety $X$. In this way we obtain
algebraic deformations over an affine variety.
These deformations are non-local, and will be the objects of our study in 
\refS{elliptic}.
\end{example}

Let $A'$ be another commutative algebra over $\K$ with a fixed
augmentation $\epsilon': A' \to \K$, and let $\phi: A \to A'$ be an
algebra homomorphism with $\phi(1) = 1$ and $\epsilon'\circ\phi =
\epsilon$. If a deformation $\lambda$ of $\L$ with base
$(A,\text{Ker\,}\epsilon = m)$ is given, then the {\it push-out}
$\lambda'=\phi_*\lambda$ is the deformation of $\L$ with base $(A',
\text{Ker\,}\epsilon' = m')$, which is the Lie algebra structure
$$
{[a_1' \textstyle{\tensor}_A (a_1 \textstyle{\tensor} l_1),
a_2' \textstyle{\tensor}_A (a_2 \textstyle{\tensor} l_2)]}_{\lambda'}
 := a'_1a'_2
\textstyle{\tensor}_A {[a_1 \textstyle{\tensor}l_1,a_2 \textstyle{\tensor}
l_2]}_{\lambda},
$$
($a'_1, a'_2 \in A', a_1, a_2 \in A, l_1, l_2 \in \L$)
on $A' \tensor \L = (A' \tensor_A A) \tensor \L = A' \tensor_A (A \tensor
\L)$.  Here $A'$ is regarded as an $A$-module with the structure $aa' =
a'\phi(a)$.

\begin{remark}
For non-local algebras there exist more than one
maximal ideal, and hence in general many different augmentations $\epsilon$.
Let $\L$ be a $\K$-vector space and assume that there exists a 
Lie $A$-algebra structure $[.,.]_A$ on $A\tensor_{\K}\L$.
Given an augmentation $\epsilon:A\to\K$ with 
associated maximal ideal  $m_\epsilon=\text{Ker\,} \epsilon$,
one obtains a Lie $\K$-algebra structure $\L^\epsilon=(\L,[.,.]_\epsilon)$
 on the
vector space $\L$ by
\begin{equation}
\epsilon \tensor \text{id\,}([1 \tensor x, 1 \tensor y]_A) = 1
\tensor [x,y]_\epsilon.
\end{equation}
Comparing this with \refD{glob} we see that 
by construction the Lie $A$-algebra $A\tensor_{\K}\L$ will be a
global deformation of the Lie $\K$-algebra  $\L^\epsilon$.
On the level of structure constants the described 
construction  corresponds simply to the effect of
``reducing the structure constants of the algebra modulo $m_\epsilon$.''
In other words, for $x,y,z\in\L$ basis elements, let
the Lie $A$-algebra structure be given by
\begin{equation}
[1\otimes x,1\otimes y]_A=\sum_z C_{x,y}^z (1\otimes z),\qquad
 C_{x,y}^z\in A.
\end{equation}
Then $\L^\epsilon$ is defined via
\begin{equation}
[x,y]_\epsilon:=\sum_z  (C_{x,y}^z \bmod m_\epsilon)\, z.
\end{equation}
In general, the algebras  $\L^\epsilon$ will be different for different 
$\epsilon$.
The Lie $A$-algebra  $A\tensor_{\K}\L$  will be a deformation 
of different Lie $\K$-algebras $L^\epsilon$.
\end{remark}
\begin{example}
For a deformation of the Lie algebra $\L=\L_0$ over the affine
line, the Lie structure $\L_\alpha$ in the fiber over the point $\alpha\in\K$ 
is given by considering the augmentation corresponding to the maximal
ideal $m_\alpha=(t-\alpha)$. 
This explains the picture in the geometric interpretation of the deformation.
\end{example}
\subsection{Formal deformations}
Let $A$ be a complete local algebra over $\K$, so
$A = \overleftarrow{\lim\limits_{n\to \infty}}(A/m^n)$, where $m$ is
the maximal ideal of $A$ and  we assume that $A/m\cong\K$.
\begin{definition}
A {\it formal deformation} of $\L$ with base $A$ is a Lie algebra
structure on the completed tensor product
$A \htensor\L =
    \overleftarrow{\lim\limits_{n\to \infty}}((A/m^n)\tensor \L)$
such that
$$
\epsilon\textstyle{\htensor}
\text{id}: A \textstyle{\htensor}\L \to \K\tensor \L = \L
$$
is a Lie algebra homomorphism.
\end{definition}

\begin{example}
If $A=\K[[t]]$, then a formal deformation of $\L$ with base $A$
is the same as a formal 1-parameter deformation of $\L$ (see \cite{Ger}).
\end{example}

There is an analogous definition for equivalence of deformations
parameterized by a complete local algebra. 
\subsection{Infinitesimal and versal formal deformations}
In the following let the base of the deformation be a local $\K$-algebra
$(A,m)$ with $A/m\cong\K$.
In addition we assume that 
$\text{dim\,} (m^k/m^{k+1}) < \infty$ for all $k$.

\begin{proposition}\label{P:infuni}
(\cite{FiaFuMini})
With the assumption $\dim \H^2(\L,\L) < \infty$, there exists a universal
infinitesimal deformation $\eta_\L$ of the Lie algebra $\L$ with
base $B= \K\oplus H^2(\L,\L)'$, where 
the second summand is the dual of $\H^2(\L,\L)$ equipped with 
zero multiplication, i.e.
$$
(\alpha_1,h_1)\cdot(\alpha_2,h_2) = (\alpha_1\alpha_2,\alpha_1h_2 +
\alpha_2 h_1).
$$
\end{proposition}
This means that for any infinitesimal deformation $\lambda$ of the Lie
algebra $\L$ with finite dimensional (local) algebra base $A$ there exists a
unique homomorphism $\phi: \K\oplus H^2(\L,\L)'\to A$ such that
$\lambda$ is equivalent to the push-out $\phi_*\eta_\L$.
\medskip

Although in general it is impossible to construct a universal 
formal deformation, there is a so-called versal element.

\begin{definition}
A formal deformation $\eta$ of $\L$ parameterized by a complete local
 algebra
$B$ is called {\it versal} if for any deformation $\lambda$,
parameterized by a complete local algebra $(A,m_A)$, there is a morphism
 $f: B\to A$ such that
\smallskip

\quad 1) The push-out $f_*\eta$ is equivalent to $\lambda$.
\smallskip

\quad 2) If $A$ satisfies ${m_A}^2=0$, then $f$ is unique
(see \cite{Fiasbor,FiaFuMini}).
\end{definition}

\begin{remark}
$ $
A versal formal deformation is sometimes called miniversal.
\end{remark}

\begin{theorem}\label{T:formver}
(\cite{Fiaproc},\cite[Thm. 4.6]{FiaFuMini}) 
Let the space  $\H^2(\L,\L)$ be finite dimensional.
\newline
(a) There exists a  versal formal  deformation of $\L$.
\newline
(b) The base of the versal formal deformation is formally embedded into
$\H^2(\L,\L)$, i.e. it can be described in $\H^2(\L,\L)$ by
a finite system of formal equations.
\end{theorem}
A Lie algebra $\L$ is called
(infinitesimally, formally, or globally) 
rigid if every (infinitesimal, formal, global) 
family is equivalent to a trivial one.
Assume  $\H^2(\L,\L)<\infty$ in the following.
By \refP{infuni} 
the elements of $H^2(\L,\L)$ correspond bijectively to the equivalence
classes of infinitesimal deformations, as equivalent deformations up to
order 1 differ from each other only in a coboundary. 
Together with \refT{formver}, Part (b), follows
\begin{corollary}\label{C:formrig}
$ $
\newline
(a) $\L$ is infinitesimally rigid if and only if  $\H^2(\L,\L)=\{0\}$.
\newline
(b)  $H^2(\L,\L) = \{0\}$ implies that $\L$ is 
formally rigid.
\end{corollary}
Let us stress the fact, that $\H^2(\L,\L)=\{0\}$ does not imply that
every global deformation will be equivalent to a 
trivial one. 
Hence, $\L$ is in this case not necessarily globally rigid.
In \refS{elliptic} we will see plenty of 
nontrivial global deformations of the Witt algebra $\W$.
Hence, the Witt algebra is not globally rigid.
In the next section we will present the proof of
 $\H^2(\W,\W)=\{0\}$, which implies infinitesimal and formal
rigidity of $\W$.
\section{Formal rigidity of the Witt and Virasoro algebras}\label{S:rigid}
As we pointed out in the Introduction, in formal deformation theory 
cohomology is an important tool.

The Lie algebras considered in this paper are infinite dimensional.
Such Lie algebras possess a topology with respect to which all algebraic
operations are continuous.  In this situation, in a cochain complex it
is natural to distinguish the sub-complex formed by the continuous
cohomology of the Lie algebra (see \cite{FeFu}).

It is known (see \cite{Fiajmp}) that the Witt and the Virasoro algebra
are formally rigid in  the sense
introduced in \refS{theory}.
The statement follows from a general result of Tsujishita
\cite{Tsu}, combined with 
results of Goncharova \cite{Gon}.
The goal of this section is to explain the relation in more detail.

First we recall the result of Tsujishita.
Recall that in this  article the Witt algebra $\W$ is defined to
be the complexification of the Lie algebra of polynomial vector fields on 
$S^1$. They constitute a dense subalgebra of the algebra 
$\text{\rm Vect\,}S^1$ of all smooth
vector fields $\text{\rm Vect\,}S^1$.
The results of Tsujishita concerns the continuous cohomology of
$\text{\rm Vect\,}S^1$ with values in formal 
tensor fields.

In fact he  deals  with the  cohomology of the algebra of vector 
fields on a general smooth compact manifold, 
but we only need his result in case of the unit circle $S^1$.
To formulate his results we have to introduce the 
space $Y(S^1)$ which is defined as follows.
Let us consider the trivial principal $U(1)$-bundle $u(S^1)$, associated
with the complexification of the real tangent bundle of $S^1$, and
let $U(S^1)$ be its total space. Denote by $x(S^1)$ the trivial
principal bundle $S^1 \times S^3 \to S^1$ with structural group $SU(2)$
and base $S^1$ and let $\Omega S^3$ be the loop space of $S^3$. The
space of sections $\text{Sec\,}x(S^1)$ 
of the bundle $x(S^1)$  is the space 
$\text{Map}(S^1,S^3) = S^3 \times
\Omega S^3$. Consider $u(S^1)$ as a subbundle of $x(S^1)$. Now $Y(S^1)$ 
is the space 
\begin{equation}
Y(S^1) := \{(y,s) \in S^1 \times \text{Sec\,}x(S^1)\,|\, s(y) \in U(S^1)\}.
\end{equation}
The space $Y(S^1)$ is homeomorphic to $S^1 \times S^1 \times \Omega S^3$,
as can be seen as follows.
We note that $s(y) \in {y} \times U(S^1)$, so we can write
the section in the form $s(u) = (u,f(u))$, $u \in S^1$, where $f: S^1
\to S^3$, $f(u) \in U(1) \subset S^3$. Now let $h$ be 
$f$  right
translated by   $f(1)^{-1}$, i.e. $h(u) = f(u)f(1)^{-1}$.
Then $h$ takes $1$ to $1$ in $S^3$ and 
we get the required mapping from $Y(S^1)$ to $S^1 \times S^1 \times
\Omega S^3$. 
On the other hand, take $y \in S^1$, $z \in S^1=U(1)$ and a loop 
$h \in \Omega S^3$ such that $h: S^1 \to S^3$, $h(1)
= 1$. Then the section $s(y) = (y, h(y)[h(y)^{-1}z])$ defines
an element of $Y(S^1)$. 

\begin{theorem} (Tsujishita \cite{Tsu}, Reshetnikov \cite{Resc})
\newline
 The cohomology ring  $\H^*(\text{\rm Vect\,}S^1, C^{\infty}(S^1))$ 
is isomorphic to $\H^*(Y(S^1),\R)$.
\end{theorem}
The real (topological) cohomology ring 
$\H^*(Y(S^1),\R)$ of the space $Y(S^1)$ 
 is known to be the free skew-symmetric $\R$-algebra 
$S(t, \theta, \xi)$, where $\text{deg\,} t = \text{deg\,} \theta
= 1, \text{deg\,} \xi = 2$.
Hence $\H^*(Y(S^1),\R)\cong S(t, \theta, \xi)$
as graded algebra.
\begin{theorem} (Tsujishita \cite{Tsu})
For an arbitrary tensor $\gl(n,\R)$-module $A$ and the space
$\mathcal{A}$ of the corresponding formal tensor fields, $\H^*(\text{\rm
Vect\,}S^1,\mathcal{A})$ is isomorphic to the tensor product of the the
ring $\H^*(Y(S^1),\R)$ and \ $\text{\rm Inv}_{\gl(n,\R)}(\H^*(L_1) \tensor
A)$, where $L_1$ denotes the subalgebra of $\W$ with basis $(l_1, l_2,
...)$. 
\end{theorem}
See the book of Fuchs \cite{Fuks} concerning this form of
the theorem and for related results.

Hence, for computing the cohomology ring, we need to know
the cohomology (with trivial coefficients)
of the Lie algebra $L_1$. This is computed by Goncharova
\cite{Gon}. She computed the cohomology spaces for all Lie algebras
$L_k$ with basis $(l_k, l_{k+1}, ...)$, but we will only state the
result we need now.  We point out that her computation is carried-out
for graded cohomology.

Let $\H_{(s)}^q$ be the $s$-homogeneous part of the cohomology space
$\H^q$ where the grading is induced by the grading of $\W$, i.e. by
$\deg l_n=n$.
\begin{theorem}\label{goncha} 
(Goncharova \cite{Gon})
For $q\ge 0$, the dimension of the graded cohomology spaces is:
\begin{equation}
\dim \H^q_{(s)}(L_1)
=\begin{cases}\quad
1,&\qquad s=\frac {3q^2\pm q}{2},
\\
\quad 0,&\qquad s\ne\frac {3q^2\pm q}{2}.
\end{cases}
\end{equation}
\end{theorem}
For the manifold $S^1$, all $\gl(1,\R)$-modules 
of formal tensor fields are of the
form $C^{\infty}(S^1)d\varphi^s$ for some integer $s$, where $\varphi$ is
the angular coordinate on the circle. Using Goncharova's and
Tsujishita's result we obtain 
\begin{theorem}
For $q\ge 0$
\begin{equation}
\H^q(\text{\rm Vect\,}S^1,C^{\infty}(S^1)d\varphi^{s})
=\begin{cases}
\quad
\H^{q-r}(Y(S^1),\R),&\qquad s=\frac {3r^2\pm r}{2},
\\
\quad\{0\},&\qquad s\ne\frac {3r^2\pm r}{2}.
\end{cases}
\end{equation}
\end{theorem}
In particular,
\begin{corollary}\label{C:vector}
 In case $s=-1$, we have
$$
\H^*(\text{\rm Vect\,}S^1,\text{\rm Vect\,}S^1) = 0.
$$
Especially, $\H^2(\text{\rm Vect\,}S^1,\text{\rm Vect\,}S^1) = 0$, so
the  algebra $\text{\rm Vect\,}S^1$ is formally rigid.
\newline
Consequently, for the algebra of complexified 
vector fields $\text{\rm Vect\,}S^1\tensor \C$ we have 
\newline 
$\H^2(\text{\rm Vect\,}S^1\tensor \C,\text{\rm Vect\,}S^1\tensor \C) = 0$,
and hence $\text{\rm Vect\,}S^1\tensor \C$
is formally rigid as well.
\end{corollary}
\begin{corollary} \label{C:witt}
$ $
\newline
(a) For the Witt algebra $\W$ we have
 $\H^2(\W,\W) =0$, hence 
the Witt algebra is  formally rigid.
\newline
(b) For the Virasoro algebra $\Vir$ we have
 $\H^2(\Vir,\Vir) =0$, hence 
the Virasoro algebra is  formally rigid
\end{corollary}
\begin{proof}
The algebra $\W$ is the subalgebra of complexified polynomial vector 
fields of $\text{\rm Vect\,}S^1\tensor \C$.
By density arguments $\H^2(\W,\W) =0$ in the graded sense and the formal
rigidity follows from 
\refC{formrig}.
The algebra $\Vir$ is a one-dimensional central extension of $\W$.
Using the Serre-Hochschild spectral sequence we obtain that
$\Vir$ as a $\Vir$-module is an extension of $\W$ as a $\W$-module.
Statement (b)  then follows from the long exact cohomology sequence.
\end{proof}

\section{Krichever-Novikov algebras}
\label{S:kn}
\subsection{The algebras with their almost-grading}
Algebras of Krichever-Novikov types are generalizations of the Virasoro
algebra and all its related algebras.
In this section we only recall the definitions and facts needed here.
Let $M$ be a compact Riemann surface of genus $g$, or in terms
of algebraic geometry, a smooth projective curve over $\C$.
Let $N,K\in\N$ with $N\ge 2$ and $1\le K<N$ be numbers. Fix
$$
I=(P_1,\ldots,P_K),\quad\text{and}\quad
O=(Q_1,\ldots,Q_{N-K})
$$
disjoint  ordered tuples of  distinct points (``marked points'',
``punctures'') on the
curve. In particular, we assume $P_i\ne Q_j$ for every
pair $(i,j)$. The points in $I$ are
called the {\it in-points}, the points in $O$ the {\it out-points}.
Sometimes we consider $I$ and $O$ simply as sets and set
$A=I\cup O$ as a set.

Denote by $\L$ the Lie algebra consisting of those meromorphic sections of the 
holomorphic tangent line bundle which are holomorphic
outside of $A$, equipped with the Lie bracket $[.,.]$ of
vector fields.
Its local form is
\begin{equation}\label{E:Lbrack}
[e,f]_|=[e(z)\fpz, f(z)\fpz]:=
\left( e(z)\pfz f(z)- f(z)\pfz e(z)\right)\fpz \ .
\end{equation}
To avoid cumbersome notation we will use the same symbol for the section
and its representing function.

For the Riemann sphere ($g=0$) with quasi-global coordinate $z$,
$I=\{0\}$ and $O=\{\infty\}$, the introduced 
vector field algebra is
the Witt algebra.
We denote for short this situation as the
{\it classical situation}.

For infinite dimensional algebras and modules
and their representation theory a graded structure is usually
of importance to obtain structure results.
The Witt algebra is a graded Lie algebra.
In our more general context the algebras will almost never be graded.
But it was observed by Krichever and Novikov 
in the two-point case that a weaker
concept, an almost-graded structure, will be enough to develop an
interesting  theory of representations (Verma modules, etc.).
\begin{definition}\label{D:almgrad}
Let $\A$ be an (associative or Lie) algebra admitting a direct
decomposition as vector space $\ \A=\bigoplus_{n\in\Z} \A_n\ $.
The algebra $\A$ is called an {\it almost-graded}
algebra if (1) $\ \dim \A_n<\infty\ $ and (2)
there are constants $R$ and  $S$ with
\begin{equation}\label{E:eaga}
\A_n\cdot \A_m\quad \subseteq \bigoplus_{h=n+m+R}^{n+m+S} \A_h,
\qquad\forall n,m\in\Z\ .
\end{equation}
The elements of $\A_n$ are called {\it homogeneous  elements of degree $n$}.
\end{definition}

For the 2-point situation for 
$M$ a higher genus Riemann surface and  $I=\{P\}$, $O=\{Q\}$
with $P,Q\in M$, Krichever and Novikov \cite{KNFac}
introduced an almost-graded structure of the vector field algebras $\L$
by exhibiting a
special basis and defining their elements to be the
homogeneous elements.
In \cite{SchlDiss,Schlkn, Schleg,Schlce} its multi-point
generalization was given, again  by
exhibiting a special basis.
Essentially, this is done by fixing their order at the points in $I$ and
$O$ in a complementary way.
For  every $n\in\Z$, and $i=1,\ldots,K$ a certain element
$e_{n,p}\in\L$ is exhibited.
The $e_{n,p}$ for $p=1,\ldots,K$ are a basis of a
subspace $\L_n$ and it is shown that
$\ \L=\bigoplus_{n\in\Z}\L_n$ .
\begin{proposition}
\cite{SchlDiss,Schlce}
With respect to the above introduced grading the 
Lie algebras $\L$ 
are almost-graded.
The almost-grading depends on the splitting $A=I\cup O$.
\end{proposition}

In the following we will have an explicit description of the basis
elements
for certain genus zero and one situation.
Hence, we will not recall their general definition.
\subsection{Central extensions}
\label{SS:cext}
To obtain the equivalent of the Virasoro algebra we have to consider 
central extensions of the algebras.
Central extensions are given by elements of $\H^2(\L,\C)$.
The usual definition of the Virasoro cocycle is not coordinate
independent. We have to introduce a projective connection $R$.
\begin{definition}
Let $\ (U_\a,z_\a)_{\a\in J}\ $ be a covering of the Riemann surface
by holomorphic coordinates, with transition functions
$z_\b=f_{\b\a}(z_\a)$.
A system of local (holomorphic, meromorphic) functions 
$\ R=(R_\a(z_\a))\ $ 
is called a (holomorphic, meromorphic) projective 
connection if it transforms as
\begin{equation}\label{E:pc}
R_\b(z_\b)\cdot (f_{\beta,\alpha}')^2=R_\a(z_\a)+S(f_{\beta,\alpha}),
\qquad\text{with}\quad
S(h)=\frac {h'''}{h'}-\frac 32\left(\frac {h''}{h'}\right)^2,
\end{equation}
the Schwartzian derivative.
Here ${}'$ denotes differentiation with respect to
the coordinate $z_\a$.
\end{definition}
\noindent
Every Riemann surface admits a holomorphic projective connection $R$
\cite{Gun}.
From \refE{pc} it follows that the difference of two projective connections
will be a quadratic differential.
Hence, after fixing one projective connection all others are obtained
by adding quadratic differentials.

For the vector field algebra $\L$
the 2-cocycle 
\begin{equation}\label{E:vecg}
\g_{S,R}(e,f):=\cintl{C_S} \left(\frac 12(e'''f-ef''')
-R\cdot(e'f-ef')\right)dz\ .
\end{equation}
defines a central extension.
Here $C_S$ is a cycle separating the in-points from the out-points.
In particular, $C_S$ can be taken to be $C_S=\sum_{i=1}^KC_i$ where
the $C_i$ are deformed circles around the points in $I$.
Recall that we use the same letter for the vector field and its
local representing function.

\begin{theorem}
\cite{Schlloc,SchlDiss}\label{T:local}
\newline
(a) The cocycle class of $\gamma_{S,R}$ does not depend on the 
chosen connection $R$.
\newline
(b) The cocycle  $\gamma_{S,R}$ is cohomologically non-trivial.
\newline
(c) The cocycle   $\gamma_{S,R}$ is local, i.e. there exists an 
$M\in\Z$ such that
\begin{equation*}
\forall n,m:\quad\gamma(\L_n,\L_m)\ne 0\implies 
M\le n+m\le 0.
\end{equation*}
(d) 
Every local cocycle for $\L$ is either a coboundary or a scalar multiple 
of \refE{vecg} with $R$ a meromorphic projective connection which
is holomorphic outside $A$.
\end{theorem}

The central extension $\Lh$ can be given
via $\Lh=\C\oplus\L$ with Lie structure 
(using the notation $\widehat{e}=(0,e)$, $c=(1,0)$)
\begin{equation}
[\widehat{e},\widehat{f}]=\widehat{[e,f]}+\g_{S,R}\,c,\quad [c,\L]=0.
\end{equation}
Using the locality, by defining $\deg c:=0$, the almost-grading can be extended
to the central extension $\Lh$.

Note that  \refT{local} does not claim that there is
only one non-trivial cocycle class (which in general is not true). It only says
that there is up to multiplication with a
scalar only one 
class such that it contains cocycles which are local with respect to
the
almost-grading. Recall that the almost-grading is  given by the splitting of 
$A$ into $I\cup O$.

%
\section{The algebra for the elliptic curve case}\label{S:elliptic}

\subsection{The family of elliptic curves}
We consider the genus one case, i.e. the case of one-dimensional
complex tori or equivalently the elliptic curve case.
We have degenerations in mind. Hence it is more convenient to use
the purely algebraic picture.
Recall that the elliptic curves can be given in the projective plane
by
\begin{equation}\label{E:ellip}
Y^2Z=4X^3-g_2XZ^2-g_3Z^3,\quad g_2,g_3\in\C,\quad
\text{with } \Delta:={g_2}^3-27{g_3}^2\ne 0.
\end{equation} 
The condition $\Delta\ne 0$ assures that the curve will be
nonsingular.
Instead of \refE{ellip} we can use the description
\begin{equation}\label{E:ellipe}
Y^2Z=4(X-e_1Z)(X-e_2Z)(X-e_3Z)
\end{equation}
with 
\begin{equation}\label{E:edef}
e_1+e_2+e_3=0,
\quad\text{and}\quad
\Delta=16(e_1-e_2)^2
(e_1-e_3)^2
(e_2-e_3)^2\ne 0.
\end{equation}
These presentations are related via
\begin{equation}
g_2=-4(e_1e_2+e_1e_3+e_2e_3),\quad
g_3=4(e_1e_2e_3).
\end{equation}
The elliptic modular parameter
classifying the elliptic curves up to isomorphy is given as 
\begin{equation}
j=1728\;\frac {g_2^3}{\Delta}
\end{equation}
We set 
\begin{equation}
B:=\{(e_1,e_2,e_3)\in \C^3\mid e_1+e_2+e_3=0,\quad
e_i\ne e_j\ \text{for}\ i\ne j\}.
\end{equation}
In the product $B\times\P^2$ we consider the family of elliptic
curves $\E$ over $B$ defined via \refE{ellipe}.
The family can be extended to 
\begin{equation}
\Bh:=\{e_1,e_2,e_3)\in \C^3\mid e_1+e_2+e_3=0
\}.
\end{equation}
The fibers above $\Bh\setminus B$ are singular cubic curves.
Resolving the one linear relation in $\Bh$ via
$e_3=-(e_1+e_2)$ we obtain a family over $\C^2$.

Consider the complex lines in $\C^2$
\begin{equation}
D_s:=\{(e_1,e_2)\in\C^2\mid 
e_2=s\cdot e_1\},\quad  s\in \C,
\qquad
D_\infty:=\{(0,e_2)\in\C^2\}.
\end{equation} 
Set also
\begin{equation}
D_s^*=D_s\setminus\{(0,0)\} 
\end{equation}
for the punctured line.
Now 
\begin{equation}
B\cong\C^2\setminus (D_1\cup D_{-1/2}\cup D_{-2}).
\end{equation}
Note that above $D_1^*$ we have $e_1=e_2\ne e_3$,
 above $D_{-1/2}^*$ we have $e_2=e_3\ne e_1$,
and  above $D_{-2}^*$ we have $e_1=e_3\ne e_2$.
In all these cases we obtain the  nodal cubic.
The nodal cubic $E_N$ can be given as
\begin{equation}
Y^2Z=4(X-eZ)^2(X+2eZ)
\end{equation}
where $e$ denotes the value of the coinciding $e_i=e_j$  ($-2e$ is then 
necessarily the remaining one).
The singular point is the point $(e: 0:1)$. It is a node.
It is up to isomorphy the only singular cubic which is stable in the
sense of Mumford-Deligne.

Above the unique common intersection point $(0,0)$ 
of all $D_s$ there is the
cuspidal cubic $E_C$
\begin{equation} 
Y^2Z=4X^3.
\end{equation}
The singular point is $(0:0:1)$. The curve is not stable
in the sense of Mumford-Deligne.
In both cases the complex projective line is the desingularisation,

In all cases (non-singular or singular) the point  
$\infty=(0:1:0)$ lies on the  
curves. It is the only intersection with the line at
infinity,  and is a non-singular point.
In  passing  to an affine chart in the following we will loose nothing.

For the  curves above the points in $D_s^*$ we calculate $e_2=s e_1$ and 
$e_3=-(1+s)e_1$ (resp. $e_3=-e_2$ if $s=\infty$).
Due to the homogeneity,  the modular parameter $j$ 
for the  curves above $D_s^*$ 
will be constant along the line.
In particular, the curves in the family lying above $D_s^*$ will
be isomorphic.
For completeness let us write down
\begin{equation}
j(s)=1728\;\frac {4(1+s+s^2)^3}{(1-s)^2(2+s)^2(1+2s)^2},
\quad j(\infty)=1728.
\end{equation}
\subsection{The family of vector field algebras}
We have to introduce the points where poles are allowed.
For our purpose it is enough to consider two marked points.
More marked points are considered in \cite{SchlDeg,RDS}.
We will always put one marking to $\infty=(0:1:0)$ and the 
other one to the point with the affine coordinate $(e_1,0)$.
These markings define two sections of the family $\E$ over $\Bh\cong
\C^2$.
With respect to the group structure on the elliptic curve
given by $\infty$ as the neutral element (the first marking) 
the second marking chooses a two-torsion point.
All other choices of two-torsion points will yield isomorphic 
situations.

In \cite{SchlDeg} for this situation (and for a three-point situation)
a basis of the Krichever-Novikov type vector field algebras were given.
\begin{theorem}\label{T:elliptic}
For any elliptic curve $E_{(e_1,e_2)}$ 
over $(e_1,e_2)\in \C^2\setminus (\Dex)$ 
the Lie algebra $\L^{(e_1,e_2)}$ of vector fields on $E_{(e_1,e_2)}$ 
has a basis $\{V_n$, $n\in\Z\}$ such that the Lie algebra structure is
given as
\begin{equation}\label{E:structell}
[V_n,V_m]=
\begin{cases}
(m-n)V_{n+m},&n,m \ \text{odd},
\\
(m-n)\big(V_{n+m}+3e_1V_{n+m-2}
\\
\qquad +\eee    V_{n+m-4}\big),&n,m \ \text{even},
\\
(m-n)V_{n+m}+(m-n-1)3e_1V_{n+m-2}
\\
\qquad +(m-n-2)\eee V_{n+m-4},&n\
\text{odd},\ m\ \text{even}.
\end{cases} 
\end{equation}
By defining $\deg(V_n):=n$, we obtain an almost-grading.
\end{theorem}
\begin{proof}
This is proved in \cite[Prop.3,Prop.4]{SchlDeg}.
Our generators are
\begin{equation}\label{E:evec}
V_{2k+1}:=(X-e_1)^{k}Y\ddX,\quad
V_{2k}:=1/2f(X)(X-e_1)^{k-2}\ddX,
\end{equation}
with $f(X)=4(X-e_1)(X-e_2)(X-e_3)$.
Note that here $V_n$ is the $V_{n-1}$ given in \cite{SchlDeg}.
\end{proof}
The algebras of \refT{elliptic}  defined with the
structure \refE{structell} make sense also for the points
$(e_1,e_2)\in D_1\cup D_{-1/2}\cup D_{-2}$.
Altogether this defines a two-dimensional family of Lie algebras
parameterized over $\C^2$.
In particular, note that we obtain for $(e_1,e_2)=0$ the Witt algebra.

Let us remark that this two-dimensional family of geometric origin
can 
also be written just with the symbols   $p$ and $q$ instead of 
$3 e_1$ and $\eee$.
In this form it was 
algebraically found by Deck \cite{Deck}, (see also Ruffing, Deck and
Schlichenmaier \cite{RDS})
as a two-dimensional family of Lie algebra.
 Guerrini \cite{Guerdef,Guerrid}  related it later 
(again in a purely algebraic manner) to  deformations of the Witt algebra 
over certain spaces of polynomials.
Due to its geometric interpretation we prefer to use the
parameterization
\refE{structell}.
Further higher-dimensional families of geometric origins 
can be obtained 
if we consider the multi-point situation 
for the elliptic curve and degenerate the curve to the cuspidal cubic
and let the marked points (beside the point at $\infty$) move to
the singularity.
But no new effects will appear.

We consider now the family of algebras  
obtained by taking  as base variety the line $D_s$ (for any $s$).
First consider $s\ne\infty$. We calculate 
$\eee=e_1^2(1-s)(2+s)$ and can rewrite for these curves 
\refE{structell} as
\begin{equation}\label{E:structs}
[V_n,V_m]=
\begin{cases}
(m-n)V_{n+m},&n,m \ \text{odd},
\\
(m-n)\big(V_{n+m}+3e_1V_{n+m-2}
\\
\qquad +e_1^2(1-s)(2+s)     V_{n+m-4}\big),&n,m \ \text{even},
\\
(m-n)V_{n+m}+(m-n-1)3e_1V_{n+m-2}
\\
\qquad
+(m-n-2)e_1^2(1-s)(2+s)    V_{n+m-4},&n\
\text{odd},\ m\ \text{even}.
\end{cases} 
\end{equation}
For $D_\infty$ we have $e_3=-e_2$ and $e_1=0$ and obtain
\begin{equation}\label{E:structi}
[V_n,V_m]=
\begin{cases}
(m-n)V_{n+m},&n,m \ \text{odd},
\\
(m-n)\big(V_{n+m}-e_2^2    V_{n+m-4}\big),&n,m \ \text{even},
\\
(m-n)V_{n+m}-
(m-n-2)e_2^2 V_{n+m-4},&n\
\text{odd},\ m\ \text{even}.
\end{cases} 
\end{equation}
If we take  $V_n^*=(\sqrt{e_1})^{-n}V_n$ (for $s\ne\infty$) as generators 
 we obtain for $e_1\ne 0$
always the algebra with $e_1=1$ in our structure equations.
For $s=\infty$ a rescaling with $(\sqrt{e_2})^{-n}V_n$
will do  the same (for $e_2\ne 0$).

Hence we see  that for fixed $s$ in all cases the algebras will be isomorphic 
above every point in $D_s$ as long as
we are not above $(0,0)$.

\begin{proposition}\label{P:noneq}
For $(e_1,e_2)\ne (0,0)$ the algebras $\L^{(e_1,e_2)}$ are not
isomorphic to the Witt algebra.
\end{proposition}
\begin{proof}
Assume that we have a Lie isomorphism $\Phi:\W=\L^{(0,0)}
\to \L^{(e_1,e_2)}$. Denote the generators of the Witt algebra by
$\{l_n$, $n\in\Z\}$.
In particular, we have $[l_0,l_n]=nl_n$ for every $n$.
We assign to every $l_n$ numbers $m(n)\le M(n)$ such that
$\Phi(l_n)=\sum_{k=m(n)}^{M(n)}\a_k(n)V_k$ with
$\a_{m(n)}(n),\a_{M(n)}(n)\ne 0$.
From the relation in the Witt algebra we obtain
\begin{equation*}
 [\Phi(l_0),\Phi(l_n)]=
\sum_{k=m(0)}^{M(0)}\sum_{l=m(n)}^{M(n)}\a_k(0)\a_l(n)[V_k,V_l]
=n\cdot \sum_{l=m(n)}^{M(n)}\a_l(n)V_l.
\end{equation*}
We can choose $n$ in such a way that the structure constants 
in the expression of $[V_k,V_l]$ at the
boundary terms will not vanish. Using the almost-graded structure
we obtain $M(0)+M(n)=M(n)$ which implies $M(0)=0$,
and $m(0)+m(n)-4=m(n)$ or $m(0)+m(n)-2=m(n)$ (for $s=1$ or 
$s=-2$) which implies 
$2\le m(0)\le M(0)=0$ which is a contradiction.
\end{proof} 
It is necessary to stress the fact, that in our approach the 
elements of the algebras are only finite linear combinations
of the basis elements $V_n$.

In particular, we obtain a family of algebras over the base  $D_s$,
which is always the affine line. In this family
the algebra over the point $(0,0)$ is the Witt algebra and
the isomorphy type above all other points will be the same 
but different from the special element, the Witt algebra.
This is a phenomena also appearing in algebraic geometry.
There it is related to non-stable singular curves (which is 
for genus one only the
cuspidal cubic).
Note that it is necessary to consider the two-dimensional family 
introduced above to ``see the full behaviour'' of the cuspidal
cubic $E_C$.

Let us collect the facts:
\begin{theorem}\label{T:families}
For every $s\in\C\cup\{\infty\}$ the families of Lie algebras defined
via the structure equations \refE{structs} for $s\ne\infty$ and
\refE{structi} for $s=\infty$ define  global deformations 
$\W^{(s)}_t$ of the
Witt algebra $\W$ over the affine line $\C[t]$. Here $t$ corresponds to
the parameter $e_1$ and $e_2$ respectively.
The Lie algebra above $t=0$ corresponds always to the Witt algebra,
the algebras above $t\ne 0$ belong (if $s$ is fixed) to the same
isomorphy type, but are not isomorphic to the Witt algebra.
\end{theorem}
 
If we denote by $g(s):=(1-s)(2+s)$ the polynomial
appearing in the structure equations \refE{structs}, we see
that the algebras over $D_s$ will be isomorphic
to the algebras over $D_t$ if $g(s)=g(t)$.
This is the case if and only if $t=-1-s$.
Under this map the lines $D_\infty$ and $D_{-1/2}$ remain fixed.
Geometrically this corresponds to interchanging the role of $e_2$ and
$e_3$.

\subsection{The degenerations and the three-point algebras for genus zero}
Next we want to identify the algebras corresponding to the 
singular cubic situation.
We have three different possibilities:
\newline
(I) All three $e_1,e_2$ and $e_3$ come together. This implies
necessarily
that $e_1=e_2=e_3=0$. We obtain the cuspidal cubic.
The pole at $(e_1,0)$ moves to the singular point $(0,0)$.
This appears if we approach in our two-dimensional 
family the point $(0,0)$.
\newline
(II) If 2 but not 3 of the $e_i$ come together 
we obtain the nodal cubic and we have to distinguish
2 subcases with respect to the marked point:
\newline
(IIa) $e_1\ne e_2=e_3$, then the point of a possible pole will remain
non-singular.
This appears if we approach a point of $D^*_{-1/2}$.
\newline
(IIb) Either $e_1=e_2\ne e_3$ or $e_1=e_3\ne e_2$, then the singular point
(the node) will become a possible pole.
This situation occurs if we approach points from $ D^*_{1}\cup D^*_{-2}$.
In the cases (IIa)  and (IIb) 
we obtain the algebras by specializing the value of $s$ in
\refE{structs}.

We want to identify  these exceptional algebras  
above $D_s$ for $s=1,-1/2$ and $-2$.

First, clearly above $(0,0)$ there is always the Witt algebra
corresponding to meromorphic vector fields on the complex line
holomorphic outside $\{0\}$ and  $\{\infty\}$.
This corresponds to situation (I)

Next  we consider the geometric situation $M=\P^1(\C)$, $I=\{\a,-\a\}$ and 
$O=\{\infty\}$, $\alpha\ne 0$.
As shown in \cite{SchlDeg} a basis 
of the corresponding Krichever-Novikov algebra is given by
\begin{equation}
V_{2k}:=z(z-\a)^k(z+\a)^k\fpz,
\qquad
V_{2k+1}:=(z-\a)^{k+1}(z+\a)^{k+1}\fpz, \qquad k\in\Z.
\end{equation}
Here $z$ is the quasi-global coordinate on $\P^1(\C)$. The grading is
given by $\deg(V_n):=n$.
One calculates
\begin{equation}\label{E:structwo}
[V_n,V_m]=
\begin{cases}
(m-n)V_{n+m},&n,m \ \text{odd},
\\
(m-n)\big(V_{n+m}+\alpha^2V_{n+m-2})
,&n,m \ \text{even},
\\
(m-n)V_{n+m}+(m-n-1)\alpha^2V_{n+m-2}
,&n\
\text{odd},\ m\ \text{even}.
\end{cases} 
\end{equation}
If we set  $\a=\sqrt{e_1}$
we get exactly the structure 
for the algebras obtained in the degeneration (IIb).
Hence,
\begin{proposition}\label{P:nodeg}
The algebras $\L^\lambda$ for $\lambda\in D_{1}^*\cup D_{-2}^*$ are
isomorphic to the algebra of meromorphic 
vector fields on $\P^1$ which are holomorphic outside
$\{\infty,\a,-\a\}$.
\end{proposition}  
Finally,  we consider the subalgebra of the Witt algebra defined by 
the basis elements
\begin{equation}
\begin{aligned}
V_{2k}&=z^{2k-3}(z^2-\a^2)^2\fpz=l_{2k}-2\a^2l_{2k-2}+\a^4l_{2k-4},
\\
V_{2k+1}&=z^{2k}(z^2-\a^2)\fpz=l_{2k+1}-\a^2l_{2k-1}.
\end{aligned}
\end{equation}
One calculates \begin{equation}
\label{E:structone}
[V_n,V_m]=
\begin{cases}
(m-n)V_{n+m},&n,m \ \text{odd},
\\
(m-n)\big(V_{n+m}-2\a V_{n+m-2}+
\a^2    V_{n+m-4}\big),&n,m \ \text{even},
\\
(m-n)V_{n+m}+(m-n-1)(-2\a)V_{n+m-2}
\\
\qquad\quad
+(m-n-2)\a^2 V_{n+m-4},&n\
\text{odd},\ m\ \text{even}.
\end{cases} 
\end{equation}
This is the algebra obtained by the degeneration (IIa)
if we set $\a=\mathrm{i}\sqrt{3e_1/2}$.
Hence,
\begin{proposition}
The algebras $\L^\lambda$ for $\lambda\in D_{-1/2}^*$ are isomorphic 
to the subalgebra of the Witt algebra generated by the above basis 
elements.
\end{proposition}
This subalgebra can be described as the subalgebra of meromorphic 
vector fields vanishing at $\a$ and $-\a$, with possible poles at 
$0$ and $\infty$ and such that in the representation of 
$V(z)=f(z)(z^2-\a^2)\fpz$ the function $f$ fulfills $f(\a)=f(-\a)$.

Clearly, as explained above, as long as $\a\ne 0$, by  rescaling,
the case $\a=1$ can be obtained. Hence for $\a\ne 0$
the algebras are all isomorphic.

If we choose a line $E$ in $\C^2$ not passing through the origin, 
by restricting our two-dimensional family to those algebras above $E$ 
we obtain a family of algebras over an affine line.
A generic line will meet $D_1,D_{-2}$ and $D_{-1/2}$.
In this way we obtain global deformations of these special algebras.
\subsection{Geometric interpretation of the deformation results}
The identification of the 
algebras obtained in the last subsection is not a pure coincidence.
There is a geometric scheme behind, which was elaborated in 
\cite{SchlDeg}.
To put the results in the right context we want to indicate the
relation.
In both cases of the singular cubic the desingularisation 
(which will be also the normalization) will be
the projective line. 
The vector fields given in \refE{evec}
make sense also for the degenerate cases.
Vector fields on the singular cubic will
correspond to vector fields on the normalization which have
at the points lying above the singular points an additional zero.

In case of the cuspidal degeneration the possible pole moves to the
singular point. Hence we will obtain the full Witt algebra.
In the case of the nodal cubics we have to distinguish the  two cases.
If $(e_1,0)$ will not be the singular point one obtains the subalgebra of
the Witt algebra consisting of vector fields which have
 a zero at $\a$ and $-\a$ (where $\a$ 
is the point lying above the singular point) 
and fulfill the additional condition on $f$ (see above).
If $(e_1,0)$ becomes a singular point, a pole at $(e_1,0)$ will
produce poles at the two points $\a$ and $-\a$ lying above  $(e_1,0)$.
Hence we end up with the 3-point algebra for genus zero.
\subsection{Cohomology classes of the deformations}
Let $\W_t$ be a one-parameter deformation of the Witt algebra $\W$
with Lie structure
\begin{equation}
[x,y]_t=[x,y]+t^k\omega_0(x,y)+t^{k+1}\omega_{1}(x,y)+\cdots,
\end{equation}
such that $\omega=\omega_0$ is a non-vanishing bilinear form.
The form $\omega$ will be a 2-cocycle in $C^2(\W,\W)$.
The element $[\omega]\in\H^2(\W,\W)$ will be the 
cohomology class characterizing the infinitesimal family.
Recall that a class $\omega$ is per definition a coboundary 
if
\begin{equation}
\omega(x,y)=(d_1\Phi)(x,y):=
\Phi([x,y])-[\Phi(x),y]-[x,\Phi(y)]
\end{equation}
for a linear map $\Phi:\W\to\W$.
For the global  deformation families  $\W^{(s)}_t$
appearing in \refT{families} we obtain with respect to 
the parameterization by  $e_1$ and $e_2$ respectively, as first
nontrivial contribution the following two cocycles.
\begin{equation}
\begin{aligned}
\omega(l_n,l_m)&=
\begin{cases}
0,&n,m\  \text{odd},
\\  (m-n)3\,l_{n+m-2},& n,m\ \text{even}
\\  (m-n-1)3\,l_{n+m-2},& n \ \text{odd}\ ,m\ \text{even},
\end{cases}
\\
\text{and}
\\
\omega(l_n,l_m)&=
\begin{cases}
0,&n,m\  \text{odd},
\\  -(m-n)l_{n+m-4},& n,m\ \text{even},
\\  -(m-n-1)l_{n+m-4},& n \ \text{odd}\ ,m\ \text{even}.
\end{cases}
\end{aligned}
\end{equation}

From \refS{rigid} we know that the 
Witt algebra is infinitesimally and formally rigid. 
Hence, the 
 cohomology classes  $[\omega]$ in 
$\H^2(\W,\W)$
must vanish. 
For illustration we will verify this directly.
By a suitable ansatz one easily finds
that $\omega=d_1\Phi$ for 
\begin{equation}
\Phi(l_n)=
\begin{cases}
-3\,l_{n-2},&n\text{\ even},
\\
-3/2\,l_{n-2},&n\text{\ odd,}
\end{cases}
\quad\text{resp.}
\quad
\Phi(l_n)=
\begin{cases}
l_{n-4},&n\text{\ even},
\\
1/2\,l_{n-4},&n\text{\ odd.}
\end{cases}
\end{equation}
From the formal rigidity of $\W$ we can conclude that the
family $\W_t^s$ considered as a formal family over
$\C[[t]]$ is equivalent to the trivial family. Hence on the formal level
there is an isomorphism $\varphi$ given by 
\begin{equation}\label{E:sumform}
\varphi_t(l_n)=V_n+\sum_{k=1}^{\infty}\alpha_kt^kV_{n-k}.
\end{equation}
Here $t$ is a formal variable.
The formal sum \refE{sumform}
will not terminate, and even if we specialize $t$ to a non-zero number,
the element $\varphi_t(l_n)$
will not live in our Krichever-Novikov algebra.
\subsection{Families of the centrally extended algebras}
In all families  considered above it is quite easy to incorporate a central
term as defined via  the local cocycle \refE{vecg}. In the genus
one case with respect to the standard flat coordinates $(z-a)$ 
the projective connection $R\equiv 0$ will do.
The difference of two projective connections will be a
quadratic differential. Hence, we obtain that any local cocycle 
is either a coboundary or can be obtained as a scalar multiple of 
\refE{vecg} with a suitable
meromorphic quadratic differential $R$ which has only poles at $A$.
The integral \refE{vecg} is written in the complex analytic picture.
But the integration over a separating cycle can be given by integration
over circles around the points where poles are allowed.
Hence, it is given as sum of residues and the 
cocycle makes perfect sense in the purely algebraic picture.
For the explicit calculation of the residue it is useful to use the
fact that for tori $T=\C/\Lambda$ with 
lattice $\Lambda=\langle 1,\tau\rangle_{\Z}$, $\mathrm{Im\,} \tau>0$, 
the complex analytic picture
is isomorphic to the algebraic elliptic curve picture via
\begin{equation}
\bar z=z\bmod \Lambda\to 
\begin{cases}
(\wp(z):\wp'(z):1),&\bar z\ne 0,
\\
(0:1:0)=\infty,&\bar z=0.
\end{cases}
\end{equation}
Here $\wp$ denotes the Weierstra\ss\ $\wp$-function.
Recall that the points where poles are allowed in the 
algebraic picture are $\infty$ and $(e_1:0:1)$.
They correspond to $\bar 0$ and $\overline{\frac 12}$
respectively in the analytic picture.
The point $\overline{\frac 12}$ is a 2-torsion point.
Replacing $\frac 12$ by any one of the other 2-torsion points
$\frac  \tau2$ and
$\frac {\tau+1}2$ respectively, yields isomorphic structures.

Let $\W_t^{(s)}$ be any of the above considered families
parameterised over $D_s^*\cong \C\setminus\{0\}$
with parameter $t$ such that $t=0$ corresponds to $(0,0)$.
Then
$\widehat{\W_t^{(s)}}=\C\oplus\W_t^{(s)}$ 
with $\widehat x=(0,x)$, $c=(1,0)$ with
\begin{equation}
[\widehat{x}, \widehat{y}]
=\widehat{[x,y]}+\widetilde{c}(t)\gamma_{S,R_t}(x,y)\cdot c,\qquad
[\widehat{x},c]=0,
\end{equation}
for $\widetilde{c}:\C\to\C$ a non-vanishing 
 algebraic function, and $R_t$ a family of quadratic
differentials varying algebraically with respect to $t$,
will define a family of centrally extended algebras.

As shown above
the non-extended algebras for fixed $s$ are mutually isomorphic. 
Recall from \refSS{cext} that choosing
a different $R_t$ gives only a different cocycle in
the same cohomology class and that $\widetilde{c}(t)$ as long as 
$\widetilde{c}(t)\ne 0$ is 
just a rescaling of the central element.
Hence, we obtain that also the centrally extended algebras are
mutually isomorphic over $D_s^*$.

The cocycle \refE{vecg} expressed as residue and calculated at
the point $\infty$ makes perfect sense also for $t=0$.
For $t=0$ it will yield the Virasoro cocycle.
In this way we obtain a 
nontrivial deformation family for the Virasoro algebra.
Clearly, to obtain examples  we might directly 
take $c\equiv 1$ and $R_t\equiv 0$.
\begin{remark}
A typical appearing  in 2-dimensional conformal field theory
of centrally extended vector field
algebras is via the Sugawara representation, i.e. by the 
modes of the energy-momentum
tensor  in representations of affine algebras 
(gauge algebras) or of algebras of $b-c$ systems.
The classical constructions extend also to the higher
genus multi-point situation,
\cite{SchlShSug,SchlDiss},
i.e. if the 2-dimensional conformal field theory is considered for
higher genus Riemann surfaces.
If one studies families of such systems varying with the moduli parameters,
corresponding to deformation of the complex structure and moving
the ``insertions points'',
one obtains in a natural way families of centrally extended algebras.
In these cases  $\widetilde{c}(t)$ and $R_t$ might vary.
In \cite{Deck} and \cite{RDS} for $b-c$-systems 
explicit formulas for the
central term  are given.
\end{remark}
\subsection{Deformations of the Lie algebra $L_1$}
Let $L_1$ be the subalgebra of the Witt algebra consisting of those vector
fields which vanish of order $\ge 2$ at $0$, i.e.
$L_1=\langle l_n\mid n\ge 1\rangle$.
It was shown by Fialowski in 
\cite{Fiarms} that this algebra is not formally rigid,
and that there are three independent formal one-parameter
deformations.
They correspond to pairwise non-equivalent deformations.
Indeed any formal one-parameter deformation of $L_1$ can be reduced
by a formal parameter change to one of these deformations
(see also \cite{FiaFusd}):
\begin{equation}\label{E:formfam}
\begin{aligned}
\  [l_n,l_m]^{(1)}_t &:= (m-n)(l_{n+m}+tl_{n+m-1});
\\
[l_n,l_m]^{(2)}_t &:= 
\begin{cases}
(m-n)l_{n+m},&n\ne 1, m\ne 1
\\
(m-1)l_{m+1}+tml_m,
&n=1, m\ne 1;
\end{cases}
\\
[l_n,l_m]^{(3)}_t &:= 
\begin{cases}
(m-n)l_{n+m},&n\ne 2, m\ne 2
\\
(m-2)l_{m+2}+tml_m,
&n=2, m\ne 2.
\end{cases}
\end{aligned}
\end{equation}
The cocycles representing the infinitesimal deformations are given
by
\begin{equation}
\begin{aligned}
\beta^{(1)}(l_n,l_m) &:= (m-n)l_{n+m-1};
\\
\beta^{(2)}(l_n,l_m)&:= 
\begin{cases}
ml_{m},&n= 1, m\ne 1
\\
0,
&n\ne 1, m\ne 1;
\end{cases}
\\
\beta^{(3)}(l_n,l_m) &:= 
\begin{cases}
ml_{m},&n= 2, m\ne 2
\\
0,
&n\ne2, m\ne 2.
\end{cases}
\end{aligned}
\end{equation}
It is shown in the above cited articles that 
the cohomology classes $[\beta^{(1)}]=[\beta^{(2)}]=0$ and
$[\beta^{(3)}]\ne 0$.
To avoid misinterpretations let us point out that these 
infinitesimal classes are not invariant under formal
equivalence of formal deformations.
Take for the first two families the Lie algebra 1-cocycles
 $\gamma^{(i)}$ with $\beta^{(i)}=d_1\gamma^{(i)}$ ($i=1,2$).
In \cite{Fiaproc}  it is shown 
that by the  formal isomorphisms $\phi_t^{(i)}(x)=x+t\gamma^{(i)}(x)$
each of these two families is  equivalent 
to a corresponding formal family for which the infinitesimal class is a 
non-vanishing scalar multiple of $[\beta^{(3)}]$.

In our geometric situation we  consider the algebra  
$W_{1,\alpha^2}:=\langle V_n\mid n\ge 1\rangle$ with the
structure equations \refE{structwo}.
If we vary $\alpha$ we obtain a family $\W_{1,\alpha^2}$.
These algebras correspond to the algebra of vector fields
on $\P^1(\C)$ which might have a pole  at the point $\infty$ and zeros
of order at least 1 at the points $\alpha$ and $-\alpha$.
By \refP{nodeg} we know that as long as $\alpha\ne 0$ they are  isomorphic
to the corresponding subalgebra of the vector field algebra $\L^\lambda$ for
$\lambda\in D_1^*\cup D_{-2}^*$.
As long as $\alpha\ne 0$ we can rescale and obtain that all these 
subalgebras belong to the same isomorphy type.
\begin{proposition}
The algebras $\W_{1,\alpha^2}$ for $\alpha\ne 0$ are not isomorphic
to the algebra $L_1$.
\end{proposition}
\begin{proof}
The proof again uses the almost-graded (respectively, graded)
structure as in the proof of \refP{noneq}.
By the absence  of $l_0$ some additional steps are
needed.
We will only sketch them.
Assume that there is an isomorphism $\phi:L_1\to\W_{1,\alpha^2}$.
From the structure \refE{structwo} we conclude that
$M(n)=nM(1)$ (notation as in the above-mentioned proof) 
with $M(1)\in\N$.
If we assume $M(1)>1$, the basis element $V_1$ will not be
in $\phi(L_1)$. So $\phi$ cannot be an isomorphism.
Hence, $M(1)=1$.
Now $\phi(l_1)=\alpha_1V_1$, $\phi(l_2)=\alpha_2V_2+\alpha_{2,2}V_1$
and in further consequence from 
$[l_1,l_2]=2l_3$ and $[l_1,l_3]=3l_4$
it follows that 
 $\phi(l_3)=\alpha_3V_3$ and $\phi(l_4)=\alpha_4V_4$.
The relations $[l_1,l_4]=3l_5$ and $[l_3,l_2]=-l_5$
in $L_1$ lead under $\phi$ to two relations in $\W_{1,\alpha^2}$
which are in contradiction.
Hence there is no such $\phi$.
\newline
Note that 
an alternative way to see the statement is to
use \refP{ident} further down.
\end{proof}
In this way we obtain a non-trivial global deformation family $\W_{1,t}$
for the algebra $L_1$.
For its 2-cocycle  we calculate
\begin{equation}
\omega(l_n,l_m)=
\begin{cases}
0,&n,m\text{ odd},
\\
(m-n)l_{n+m-2},&n,m\text{ even},
\\
(m-n-1)l_{n+m-2},&n\text{ odd},\ m\text{ even}.
\end{cases} 
\end{equation}
Again with a suitable ansatz we find  with
\begin{equation}
\Phi(l_m):=
\begin{cases}
-\frac {m+8}{6}\,l_{m-2},&m\text{ even},\  m\ge 4,
\\
-\frac {m+5}{6}\,l_{m-2},&m\text{ odd},\  m\ge 3,
\\
0,&m=1,m=2.
\end{cases} 
\end{equation}
that $\omega-d_1\Phi=\frac 13\beta^{(3)}$.

From the structure equations \refE{structwo} we can immediately verify 
\begin{lemma}\label{L:commu}
For $\alpha\ne 0$ the commutator ideal calculates to
\begin{equation}
[W_{1,\alpha^2},W_{1,\alpha^2}]=\langle V_n\mid n\ge 3\rangle,
\end{equation}
and we have 
\begin{equation}
\dim W_{1,\alpha^2}/[W_{1,\alpha^2},W_{1,\alpha^2}]=2.
\end{equation}
\end{lemma}
\begin{proposition}\label{P:ident}
The family $\W_{1,\alpha^2}$ is formally equivalent to the first
family $[.,.]_t^1$ in \refE{formfam}.
\end{proposition}
\begin{proof}
Clearly  $\W_{1,\alpha^2}$ defines a formal family. From the 
above calculation we obtain $[\omega]=1/3[\beta^{(3)}]$.
But $[\beta^{(3)}]\ne 0$, hence 
considered as formal family it is also non-trivial.
By the results about the versal family of $L_1$ it has
to be equivalent to one of the three families of \refE{formfam}.
Only the first family $L_1^{(1)}$ has 
$\dim L_1^{(1)}/[L_1^{(1)},L_1^{(1)}]=2$,  for the other two this
dimension equals one.
By \refL{commu} we obtain the claim.
\end{proof}



\end{document}